\documentclass[a4paper,draft,reqno,12pt]{amsart}

\usepackage{tikz}
\usepackage{amssymb}
\usepackage{enumerate}
\usepackage{euscript}
\usepackage[all]{xy}
\usepackage[shortcuts]{extdash}
\newtheorem{theorem}{Theorem}[section]

\newtheorem{proposition}[theorem]{Proposition}

\theoremstyle{definition}
\newtheorem{definition}[theorem]{Definition}
\newtheorem{example}[theorem]{Example}

\newtheorem{corollary}[theorem]{Corollary}

\theoremstyle{remark}
\newtheorem{remark}[theorem]{Remark}

\numberwithin{equation}{section}

\newcommand{\rn}[1]{\MakeUppercase{\romannumeral #1}}
\newcommand{\K}{\mathbb{K}}
\newcommand{\Z}{\mathbb{Z}}
\newcommand{\D}{\mathfrak{D}}
\newcommand{\Q}{\mathbb{Q}}

\renewcommand{\O}{\mathcal{O}}

\newcommand{\Zo}{\Z_{>0}}
\newcommand{\ra}{\rightarrow}
\renewcommand{\P}{\mathbb{P}}

\renewcommand{\c}{\text{{\upshape cone}}}
\newcommand{\im}{\text{{\upshape Im}}}

\newcommand{\supp}{\text{{\upshape Supp}}}
\newcommand{\lcm}{\text{{\upshape lcm}}}
\newcommand{\pol}{\text{{\upshape Pol}}}
\newcommand{\dv}{\text{{\upshape Div}}}
\newcommand{\cdv}{\text{{\upshape CaDiv}}}
\newcommand{\hm}{\text{{\upshape Hom}}}
\newcommand{\spc}{\text{{\upshape Spec}}}
\renewcommand{\ker}{\text{{\upshape Ker}}}
\newcommand{\conv}{ \text{{\upshape conv}}}
\newcommand{\pp}{ \text{{\upshape PPDiv}}}
\newcommand{\rk}{ \text{{\upshape rank}}}
\newcommand{\hyph}{-\penalty0\hskip0pt\relax}
\newcommand{\ord}{ \text{{\upshape ord}}}

\begin{document}

\title{Polyhedral divisors of affine trinomial hypersurfaces}

\author{Oleg Kruglov}
\address{Lomonosov Moscow State University, Faculty of Mechanics and Mathematics, Department
of Higher Algebra, Leninskie Gory 1, Moscow, 119991 Russia}
\email{o.kruglov@mail.ru}





\begin{abstract}
We find polyhedral divisors corresponding to the torus action of complexity one on affine trinomial hypersurfaces. Explicit computations for particular classes of such hypersurfaces including Pham\=/ Brieskorn surfaces and rational trinomial hypersurfaces are given.
\end{abstract}

\maketitle

\section{Introduction}
One of the main problems of modern algebraic geometry is to find an effective description of specific classes of algebraic varieties and to study the geometric properties of varieties in terms of this description. A classical example of the realization of such an approach is the theory of toric varieties. In this case, a normal algebraic variety with the action of an algebraic torus with an open orbit is given by a fan of rational polyhedral cones, see, for example,~\cite{CLS, F}. Many results have been obtained that relate the geometry of a toric variety to the combinatorial properties of the corresponding fan.

A natural generalization of this theory is the study of arbitrary actions of an algebraic torus. An important invariant characterizing the action is the complexity of the action, defined in 1986 by Vinberg~\cite{Vin}. In the case of an algebraic torus it is equal to the codimension of a generic orbit. Proper polyhedral divisors (pp-divisors) were introduced as a tool for description of $T$\=/ varieties i.e. normal affine varieties with torus actions of an arbitrary complexity by Altmann and Hausen~\cite{AH}. This approach has become widespread. For more details, the reader can refer to the survey~\cite{TVar}. Thus, the problem of an explicit computation of polyhedral divisors for concrete affine varieties is of interest. This description is the most effective in the case of actions of complexity one.

The theory of Cox rings associates arbitrary actions of complexity one with affine varieties given by trinomial equations. The case of a single trinomial equation leads to affine trinomial hypersurfaces. The problem of computing pp-divisors corresponding to the action of a torus on a trinomial hypersurface was previously discussed in \cite{Arzh}. Namely, Arzhantsev has found pp-divisors corresponding to the action of a torus on factorial trinomial hypersurfaces.

In this paper we find a general form of a polyhedral divisor corresponding to the canonical action of a torus of complexity one on an affine trinomial hypersurface. In the last section we give an explicit calculation of the polyhedral divisor for rational trinomial hypersurfaces. A classification of such hypersurfaces and a general form of the divisor for each of the classes are given. A two-dimensional case is considered separately, corresponding affine surfaces are widely known as Pham-Brieskorn surfaces. In the final example, we consider the case of a non-rational hypersurface; in this example a curve of genus one appears.

The author is grateful to I.~V.~Arzhantsev for the formulation of an interesting problem and a lot of useful conversations.

\section{Polyhedral divisors}
An irreducible closed codimension\=/ one subvariety of an algebraic variety $Y$ is called {\it a~prime divisor}. {\it A Weil divisor} on $Y$ is a finite formal sum over the prime divisors of~$Y$,
\[
D = \sum\limits_i n_i \cdot D_i,
\]
where $n_i \in \Z$. The union of all prime divisors $D_i$ in this sum with $n_i \not= 0$ is called {\it the support of divisor} $D$ and denoted by $\supp(D)$. The group of all Weil divisors on $Y$ is denoted by $\dv(Y)$. We can also define the group of rational Weil divisors $\dv_\Q(Y) := \Q \otimes_\Z \dv(Y)$. A~Weil divisor is {\it effective} if all the coefficients $n_i$ are non-negative. For every non\=/ zero rational function $f$ on $Y$ one can define the Weil divisor as follows: 
\[
(f):= \sum\limits_i \ord_{D_i}(f),
\]
where $\ord_{D_i}(f)$  is the number of zeroes or poles at $D_i$ (and is negative if it represents the number of poles). The divisor $ D $ obtained from the function $ f $ in the way described above, is called {\it principal}. {\it A~Cartier divisor} on a variety $ Y $ is a locally principal divisor, that is, a divisor of $ D $ such that there exists an open covering $ Y $ of affine $ U_i $ such that for any $ i $, the restriction of $ D $ to $ U_i $ coincides with the restriction on $ U_i $ of the principal divisor of some function $ f_i \in \ K (Y) $. Cartier divisors on $ Y $ form a subgroup of $ \dv (Y) $, denoted by $ \cdv (Y) $ and, respectively, $ \cdv_ \Q (Y): = \Q \otimes_\Z \cdv (Y)$. 

Divisor $ D \in \dv (Y) $ can be associated with a subspace in $ \K (Y) $  defined as
\[
\Gamma(Y, \O(D)) := \{f \in \K(Y)^\times\, \mid\,  (f) + D \ge 0\} \cup \{0\}.
\]
Consider $f \in \Gamma(Y, \O(D))$. Its zero set can be defined as $Z(f) := \supp(D + (f))$. Hence, principal open subset given by $f$ can be understood as $Y_f := Y\backslash Z(f)$. The divisor $D$ is called {\it semiample} if there is $n \in \Zo$ such that $Y$ is covered by the union of principal open subsets $Y_f$ for all $f\in \Gamma(Y, \O(nD))$. The divisor $D$ is called {\it big} if there is $n \in \Zo$ such that $Y_f$ is affine for some $f\in \Gamma(Y, \O(nD))$.

We need some basic concepts from convex geometry to completely explain definition of a polyhedral divisor. {\it Convex polyhedral cone} generated by vectors $\{v_1, \dots , v_n\}$ in a rational vector space $V$ is set of linear combinations of these vectors with non-negative coefficients. Further we call it simply a cone for brevity. The intersection of finitely many affine half-spaces in affine space $\K^n$ is called {\it a polyhedron}. For a polyhedron $ \Delta $, we define the set of relative interior points which are interior points of the polyhedron lying in the minimal affine subspace containing $\Delta$. Denote it as $\text{relint}(\Delta) $. Let $N$ be a lattice and $N_\Q := \Q \otimes_\Z N$ be a rational vector space associated with the lattice $N$. We let $M := \hm(N, \Z)$ be the dual lattice to $N$ and similarly $M_\Q := \Q \otimes_\Z M$. For every cone $\sigma$ in $N_\Q$ there is a dual cone $\sigma^\vee \subseteq M_\Q$ defined as
\[
\sigma^\vee := \{u \in M_\Q\, |\, \forall v \in \sigma \: \: \langle u, v\rangle \ge 0\}.
\]
An important characteristic of a polyhedron $\Delta$ is the {\it recession cone}. It is the cone consisting of rays from a fixed point of the polyhedron which are fully contained in $\sigma$. It is an easy exercise to show that a recession cone does not depend on a point choice and it is a cone. All polyhedra with a fixed recession cone $\sigma$ form a commutative monoid with the Minkowski sum as an operation. A Grothendieck group of this monoid is denoted by $\pol_\sigma (N_\Q)$.
\begin{definition}
{\it A polyhedral Weil divisor} associated with the cone $\sigma$ (or {\it {$\sigma$\=/ polyhedral} divisor}) is an element of the group
\[
\dv_\Q(Y, \sigma) := \pol_\sigma(N_\Q) \otimes_\Z \dv(X).
\]
\end{definition}
Consider $u \in \sigma^\vee$. It can be interpreted as a linear function on $N_\Q$, $v \mapsto \langle u, v\rangle$. It is well known that the minimum of this function on every polyhedron $\Delta \in \pol_\sigma (N_\Q)$ is finite. This fact allows us to define the following linear map for every $u \in \sigma^\vee$:
\begin{gather*}
\dv_\Q(Y, \sigma) \ra \dv_\Q(Y),\\
\D = \sum\limits_i \Delta_i \cdot D_i \mapsto \D(u) := \sum\limits_i \min\limits_{v \in \Delta_i}\langle u, v\rangle \cdot D_i.
\end{gather*}
We are now ready to define proper polyhedral divisors.
\begin{definition}
A polyhedral divisor $\D \in \dv_\Q(Y, \sigma)$ is called {\it a proper polyhedral divisor} if it can be represented as
\[
\D= \sum\limits_i \Delta_i \cdot D_i,
\]
where $\Delta_i \in \pol^+_\sigma(N_\Q)$ and $\D(u) \in \cdv_\Q(Y)$; moreover, $\D(u)$ is semiample for all $u \in \sigma^\vee$ and big for all $u \in \text{relint}(\sigma^\vee)$.
\end{definition}
Proper polyhedral divisors form semigroup $\text{PPDiv}_\Q(Y, \sigma) \subseteq \dv_\Q(Y, \sigma)$. If $Y$ is semiprojective, i.e. it is projective over some affine variety, it is possible to establish the correspondence between proper $\sigma$\=/ polyhedral divisors on $Y$ and $T$\=/ varieties.
\begin{theorem}~\cite{AH}
To any proper $\sigma$\=/ polyhedral divisor $\D$ on a semiprojective variety $Y$ one
can associate a normal affine $T$\=/ variety of dimension $\rk(M) + \dim(Y)$ given by $X[Y,D] =
\spc(A[Y, \D])$, where
\[
A[Y, \D] = \bigoplus\limits_{m \in \sigma^\vee \cap M} \Gamma\left(Y, \O(\D(m))\right).
\]
Conversely, any normal affine $T$\=/ variety is isomorphic to $X[Y,\D]$ for some semiprojective
variety $Y$ and some proper $\sigma$\=/ polyhedral divisor $\D$ on $Y$.
\end{theorem}
So, to determine an effective torus action on a normal affine variety we need to find a semiprojective variety $Y$ and a proper $\sigma$\=/ polyhedral divisor on it.
We call $Y$ the base variety.
\section{Main results}
Let $\K$ be algebraically closed field and $\text{char}(\K) = 0$. Fix positive integer numbers $n_0, n_1, n_2$ and $n = n_0 + n_1 + n_2$. For each $n_i$ consider a tuple $l_i \in \Z_{\ge 0}^{n_i}$ and define
\[
T_i^{l_i}:=T_{i1}^{l_{i1}} \dots T_{in_i}^{l_{in_i}} \in \K[T_{ij}; i = 0,1,2, j = 0, \dots n_i].
\]
The hypersurface $X$ is called {\it trinomial} if it is the zero set of the polynomial ${T_0^{l_0} + T_1^{l_1} + T_2^{l_2}}$ in the affine space $\K^n$. We suppose here that $n_i l_i > 1$ for $i = 0, 1, 2$, otherwise such a hypersurface is isomorphic to $\K^{n-1}$. Any trinomial hypersurface admits a torus action of complexity one. It means a generic orbit of such action has codimension one. One can define this action as follows: consider a restriction of the natural action of an $n$\=/ dimensional torus $\bar T$ on $\K^n$. Let $\bar N$ be a one-parameter subgroup lattice of the torus $\bar T$ and $\{e_1, \dots , e_n\}$ be a basis of $\bar N$. The torus $T$ is a stabilizer of the ideal $I(X) \subseteq \K[X]$. This torus can be descripted by an integer \mbox{$(2 \times n)$\=/ matrix} 
\[
L = \begin{pmatrix} -l_0 & l_1 & 0 \\ -l_0 & 0 & l_2 \end{pmatrix}.
\]
It determines the linear map $L: \bar N \ra \Z^2$. It is easy to see that $T$ is a torus with a one-parameter subgroup lattice $N = \ker(L)$. So, an action of the torus $T$ arises on~$X$. Let $F: N \rightarrow \bar N$ be the inclusion mapping. We need to choose a map $S: \bar N \rightarrow N$, such that $S \circ F = \text{id}_N$.
Note $S$ is not uniquely defined.
There is $M$\=/ grading on $\K[X]$ corresponding to this action, where $M = \hm(N, \Z)$ is a character lattice of the torus $T$. Let us illustrate the described construction with an example.
\begin{example}\label{start}
Consider a trinomial hypersurface $X \subseteq \K^4$ which is defined by the equation $T_{01}^3 + T_{11}^5 + T_{21}T_{22} =0$. In this case, the matrix $ L $ has the form
\[
L= \begin{pmatrix}-3 & 5 & 0 & 0\\-3 & 0 & 1 & 1\end{pmatrix}.
\]
The kernel $ N $ of this mapping is generated by vectors $(5, 3, 0 , 15)$ and $ (0, 0, 1, -1) $. So, the action of a two-dimensional torus on $X$ looks like
\[
(t_1, t_2) \cdot (T_{01}, T_{11}, T_{21}, T_{22}) = (t_1^5T_{01}, t_1^3T_{11}, t_2T_{21}, t_1^{15}t_2^{-1}T_{22}),
\]
and $M$-grading of $\K[X]$ given by this action is determined by
\[
\deg(T_{01}) = (5, 0), \: \deg(T_{02}) = (3, 0), \: \deg(T_{11}) = (0,1), \: \deg(T_{21}) = (15,-1).
\]
\end{example}
We introduce the notation that will be needed to formulate the main result of this paper.
\begin{gather*}
d_i := \gcd(l_{i0}, \dots , l_{in_i}), \, d := \gcd(d_0, d_1, d_2),\\
d_{ij} := \gcd\left(\frac{d_i}{d}, \frac{d_j}{d}\right), \,\tilde d := dd_{01}d_{02}d_{12}.
\end{gather*}

%
%
\begin{theorem}\label{main}
Let $X \subseteq \K^n$ be a trinomial hypersurface. Then the pp-divisor corresponding to a natural action of an $(n - 2)$-dimensional torus on $X$ is equal~to
\[
\D = \Delta_0 \cdot D_0 + \Delta_1 \cdot D_1 + \Delta_2 \cdot D_2,
\]
where $D_i, \, i = 0, \, 1, \, 2$ are divisors on a plane curve $Y$ in a weighted projective space $\P(d_{12}, d_{02}, d_{01})$. The curve is a zero set of the homogeneous polynomial 
\[
w_0^{\tilde d / d_{12}} + w_1^{\tilde d / d_{02}} + w_2^{\tilde d / d_{01}} = 0,
\] 
and $D_i = Y \cap \{[w_0:w_1:w_2] \in \P(d_{12}, d_{02}, d_{01}) \mid w_i = 0\}$. Polyhedra have the recession cone $S(\Q^n_{\ge 0} \cap N_\Q)$ and their sets of vertices are
\begin{gather*}
V(\Delta_0) = S\left(\frac{dd_{01}d_{02}}{l_{0j}}e_k\right), \, j = k, \,  k = 1, \dots , n_0,\\
V(\Delta_1) = S\left(\frac{dd_{01}d_{12}}{l_{1j}}e_k\right), \, j = k - n_0, \, k = n_0 + 1,\, \dots, \,n_0 + n_1,\\
V(\Delta_2) = S\left(\frac{dd_{02}d_{12}}{l_{2j}}e_k\right), \, j =  k - n_0 - n_1, \, k = n_0 + n_1 + 1, \,\dots,\, n.
\end{gather*}
\end{theorem}
\section{Proof of Theorem~\ref{main}}
Let us recall a method of calculating pp-divisors for normal affine $T$\=/ varieties \mbox{from~\cite[Section~11]{AH}}. Consider an action $T \times X \rightarrow X$ as a restriction of the diagonal action of $T$ on $\K^n$. Thus, $T$ is a subtorus of the torus $\bar T$ of all invertible diagonal matrices. The inclusion $T \subseteq \bar T$ corresponds to the inclusion $F: N \subseteq \bar N$ of lattices of one-parameter subgroups. We obtain a split exact sequence
\[
0 \longrightarrow N \stackrel{F}{\longrightarrow} \bar N \stackrel{P}{\longrightarrow} N' \longrightarrow 0,
\]
where
$N' = \bar N / N$,
and
$P$
is a canonical projection. Let $\Sigma$ be the coarsest common refinement of cones $P(\sigma)$, where $\sigma$ runs through all faces of the cone $\Q^n_{\ge 0} \subseteq \bar N_{\Q}$. The toric variety $W$ defined by $\Sigma$ is a base variety and $O$ is an open orbit of the torus action. Let us denote the set of all one-dimensional cones in $\Sigma$ as $R$ and $D_\rho$ is the prime divisor corresponding to $\rho$. For $\rho \in R$ consider $v_\rho$ which is the generator of the semigroup $\rho \cap M$. Define the polyhedron
\[
\Delta_\rho := S\left(\Q^n_{\ge 0} \cap P^{-1}(v_\rho)\right) \subseteq N_\Q.
\]
Then the proper polyhedral divisor $\D_{toric}$ for action of the torus $T$ on $\K^n$ looks as follows:
\[
\D_{\text{toric}} = \sum_{\rho \in R} \Delta_\rho \cdot D_\rho.
\]
We can come now to the torus action on $X$. In this case a pp-divisor can be obtained as a restriction of the divisor $\D_{toric}$ on a normalization of an image of $X \cap O$ under the map~$P$.
We apply this construction to the trinomial case. The map $P$ is defined by the matrix~$L$. Note that $L$ can be non-surjective, so we assume $N' = \im(L)$. Hence, $N'$ is generated by vectors
\[
v_0 = (-d_0, -d_0), \, v_1 = (d_1, 0), \, v_2 = (0, d_2).
\]
The fan $\Sigma$ is generated by vectors $v_0$, $v_1$, $v_2$ as well in sense that for every proper subset of $\{v_0, v_1, v_2\}$ there is a cone generated by this subset in $\Sigma$.
\begin{proposition}\label{latgen}
Let $N$ be the lattice generated by vectors $\{v_0, \dots, v_n\}$ such that $\sum_{i = 0}^n q_i v_i = 0, \, q_i \in \Z_{>0}$, $\gcd(q_0, \dots, q_n) = 1$. Then $v_i = a_iu_i$, $i = 0, \dots, n$ where $u_i$ is a generator of the semigroup $\Q_{\ge 0} v_i \cap N$ and
\[
a_i := \gcd(q_0, \dots, q_{i-1}, q_{i+1}, \dots, q_n).
\]
\end{proposition}
One can find the proof of Proposition~\ref{latgen} in~\cite[Lemma~1(c)]{WPS}. According to this proposition, the lattice $N'$ is generated by vectors
\[
u_0 = dd_{01}d_{02}\cdot(-1, -1), \, u_1 = dd_{01}d_{12}\cdot(1, 0), \,  u_2 = dd_{02}d_{12}\cdot (0,1).
\]
The following proposition leads us to finding the toric variety $W$.
\begin{proposition}
Let $N$ be a lattice and $\{e_1, \dots e_n\}$ be a basis of $N$. Consider vectors $\{v_0, \dots, v_n\}$ such that 
\[
v_0 = -\frac{1}{q_0}\sum_{i=1}^ne_i,\,v_i = \frac{1}{q_i}e_i, \, i = 1, \dots , n,
\]
where $q_i \in \Z_{>0}$, $i = 0, \dots, n$. Then the toric variety corresponds to the fan generated by $\{v_0, \dots v_n\}$ is the weighted projective space $\P(q_0, \dots , q_n)$.
\end{proposition}
A more detailed explanation of this fact is given in~\cite[Section~2.3]{F}. Therefore, $W$ is $\P(d_{12}, d_{02}, d_{01})$ and $\D_\text{toric}$ looks like
\[
\D_{\text toric} = \Delta_0 \cdot D_0 + \Delta_1 \cdot D_1 + \Delta_2 \cdot D_2 ,
\]
where $D_i$ is a prime divisor given by the equation $w_i = 0$. Let us describe polyhedral coefficients of $\D_\text{toric}$. They have the form
\[
\Delta_i = S\left(\Q^n_{\ge 0} \cap L^{-1}(u_i)\right), \, i = 0,\,1,\,2.
\]
Since $\ker(L)=N_\Q$, a recession cone of polyhedra $\Delta_i$ is the cone $\sigma = S(\Q^n_{\ge 0} \cap N_\Q)$. Consider images of basis vectors of $\bar N_\Q$ under the mapping $L$:
\begin{gather*}
L(e_k) = \frac{l_{0j}}{dd_{01}d_{02}}u_0, \, j = k, \, k = 1, \dots , n_0;\\
L(e_k) =  \frac{l_{1j}}{dd_{01}d_{12}}u_1, \, j = k - n_0, \, k = n_0 + 1, \:\dots, \: n_0 + n_1;\\
L(e_k) = \frac{l_{2j}}{dd_{02}d_{12}}u_2, \, j = k - n_0 - n_1, \, k = n_0 + n_1 + 1, \: \dots,\: n.
\end{gather*}
It is easy to see that
\begin{gather*} 
\Q^n_{\ge 0} \cap L^{-1}(u_0) = \left\{\left.\sum_{k=1}^{n_0}\lambda_k\frac{dd_{01}d_{02}}{l_{0j}}e_k + n \,\right| \, \lambda_k \ge 0, \,\sum_{k=1}^{n_0}\lambda_k = 1, \,n \in N_\Q\right\},\\
\Q^n_{\ge 0} \cap L^{-1}(u_1) = \left\{\left.\sum_{k = n_0 + 1}^{n_0 + n_1}\lambda_k\frac{dd_{01}d_{12}}{l_{1j}}e_k + n \,\right| \,\lambda_k \ge 0 , \,\sum_{k = n_0 + 1}^{n_0 + n_1}\lambda_k =1, \, n \in N_\Q \right\},\\
\Q^n_{\ge 0} \cap L^{-1}(u_2) = \left\{\left.\sum_{k = n_0 + n_1 + 1}^n \lambda_k\frac{dd_{02}d_{12}}{l_{2j}}e_k + n \,\right|\, \lambda_k \ge 0, \, \sum_{k = n_0 + n_1 + 1}^n \lambda_k = 1, \, n \in N_\Q \right\}.
\end{gather*}
Hence, polyhedra $\Delta_i$ have following sets of vertices:
\begin{gather*}
V(\Delta_0) = \left\{\left.S\left(\frac{dd_{01}d_{02}}{l_{0j}}e_k\right) \right| j = k, \,  k = 1, \dots , n_0\right\},\\
V(\Delta_1) = \left\{\left.S\left(\frac{dd_{01}d_{12}}{l_{1j}}e_k\right) \right| j = k - n_0, \, k = n_0 + 1,\, \dots, \,n_0 + n_1\right\},\\
V(\Delta_2) = \left\{\left.S\left(\frac{dd_{02}d_{12}}{l_{2j}}e_k\right) \right| j =  k - n_0 - n_1, \, k = n_0 + n_1 + 1, \,\dots,\, n\right\}.
\end{gather*}
So, polyhedral coefficients of the divisor $\D_\text{toric}$ is
\[
\Delta_i = S\left(\text{conv}(V(\Delta_i) +  \Q_{\ge 0}^n \cap N_\Q\right), \, i = 0, 1, 2.
\]
Now consider the map $L$ on tori $\bar T \subseteq \K^n$ and $(\K^\times)^2 \subseteq \P(d_{12}, d_{02}, d_{01})$.
It can be represented as a composition of maps of lattices
\[
\bar N \stackrel{L}{\longrightarrow} N_1 \longrightarrow N_2 \stackrel{\phi^{*}}{\longrightarrow} N',
\]
where $N_1 = \Z^2$ is the lattice containing $N'$ and generated by vectors $\{(1,0), \,(0,1)\}$. Consider $\Sigma$ as a fan in $N_1$. Its maximal cones are
\[
\c\left((1, 0), \,(0, 1)\right), \, \c\left((1, 0),\,(-1,-1)\right), \, \c\left((0,1), \,(-1,-1)\right),
\]
and this fan corresponds to $\P^2$. The mapping of related tori looks like
\begin{gather*}
(t_0,\, t_1, \, t_2) \mapsto \left[1: t_0^{-l_0}t_1^{l_1}: t_0^{-l_0}t_2^{l_2}\right] = \left[ t_0^{l_0}: t_1^{l_1}: t_2^{l_2}\right], \\
\text{where }t_i = (t_{i1}, \dots , t_{in_i}), \, t_i^{l_i} = t_{i1}^{l_{i1}}, \dots , t_{in_i}^{l_{in_i}}.
\end{gather*}
The image of $X \cap O$ under this map is given by the equation
\[
{x_0 + x_1 + x_2 = 0}.
\]
The lattice $N_2$ is generated by vectors $\{(\tilde d, 0),\, (0, \tilde d)\}$, the corresponding toric morphism~is
\[
[x_0:x_1:x_2] \mapsto [x_0^{\tilde d}:x_1^{\tilde d}:x_2^{\tilde d}],
\]
and, consequently, a closure of an image of $X\cap O$ under this map goes into the~curve
\[
x_0^{\tilde d} + x_1^{\tilde d} + x_2^{\tilde d} = 0.
\]
Finally, the map $N_2 \rightarrow N'$ corresponds to a morphism of factorization on $\P^2$ by the action of the group $N' / N_2 \simeq \mu_{d_{12}} \oplus \mu_{d_{02}} \oplus \mu_{d_{01}}$ where $\mu_k$ is the group of $n$\=/ th roots of unity. The action is defined as follows:
\begin{gather*}
(\zeta_0,\zeta_1,\zeta_2) \cdot [w_0:w_1:w_2] = [\zeta_0 w_0:\zeta_1 w_1:\zeta_2 w_2], \\
\text{where}\,(\zeta_0,\zeta_1,\zeta_2) \in \mu_{d_{12}} \oplus \mu_{d_{02}} \oplus \mu_{d_{01}}.
\end{gather*}
Let us note that the quotient isomorphism has the form
\begin{gather*}
\phi: \P_2/(\mu_{d_{12}} \oplus \mu_{d_{02}} \oplus \mu_{d_{01}}) \ra \P(d_{12}, d_{02}, d_{01}),\\
x_0^{d_{12}} \mapsto w_0, \, x_1^{d_{02}} \mapsto w_1, \, x_2^{d_{01}} \mapsto w_2.
\end{gather*}
There we have used the proposition which is proved in~\cite[Proposition~1.3.18]{CLS}
\begin{proposition}
Let $N'$ be a sublattice of a finite index in $N$ with the quotient $G = N/N'$ and
$\sigma \subseteq N'_\Q = N_\Q$ be a strongly convex cone. Then
\begin{enumerate}
\item there are natural isomorphisms
$G \simeq \hm_\Z(M'/M, \K^\times) = \ker(T' \rightarrow T)$;
\item $G$ acts on $\K[\sigma^\vee \cap M']$ with ring of invariants
$\K[\sigma^\vee \cap M']^G = \K[\sigma^\vee \cap M]$;
\item $G$ acts on $U_{\sigma, N'} := \spc\left(\K[\sigma^\vee \cap M']\right)$ and the morphism $\phi: U_{\sigma, N'} \ra U_{\sigma, N}$ is constant on $G$\=/ orbits
and induces a bijection
$U_{\sigma, N'}/G \simeq U_{\sigma, N}$.
\end{enumerate}
\end{proposition}
As a result we have got the closure of an image of $X \cap O$ under the map $L$. It is the curve given by the equation $w_0^{dd_{01}d_{02}} + w_1^{dd_{01}d_{12}} + w_2^{dd_{02}d_{12}} = 0$. Now we need to check that this curve is normal to complete the proof.

Consider $\Z$\=/ graded algebra $\K[w_0, w_1, w_2]$, the grading arises from
\[
\deg w_0 = d_{12}, \, \deg w_1 = d_{02}, \, \deg w_2 = d_{01},
\]
Let $\pi$ be the graded morphism from $\K[w_0, w_1, w_2]$ to polynomials with the natural $\Z$\=/ grading
\begin{gather*}
\pi: \K[w_0, w_1, w_2] \ra \K[x_0, x_1, x_2],\\
w_0 \mapsto x_0^{d_{12}}, \, w_1 \mapsto x_1^{d_{02}}, \, w_2 \mapsto x_2^{d_{01}}.
\end{gather*}
\begin{definition}
A straight cover of the projective curve $Y$ given by an equation $f = 0$,  $f \in \K[w_0, w_1, w_2]$ is the curve $\tilde{Y} \in \P^2$ given by the equation $\pi(f)=0$.
\end{definition}
It is a quotient by an action of a finite group. Hence, if $\tilde Y$ is normal then $Y$ is normal. A curve is normal if and only if it is smooth. So, we need to check if the projective curve $\tilde Y$ smooth. Consider the affine patch  $U_0: \,x_0 = 1$. The image of $\tilde Y$ on this patch is $1 + y_1^{\tilde d} + y_2^{\tilde d} = 0$, where $y_1 = \frac{x_1}{x_0}$, $y_2 = \frac{x_2}{x_0}$. It is obvious this affine curve is smooth for all~$\tilde{d}$.
Other two patches can be considered similarly. Thus, the curve $Y$ is smooth. It completes the proof.
\begin{remark}
Let us consider a structure of divisors $D_i$ in more detail. The inverse image of the divisor $D_i$ by the map $\pi$ looks like
\[
\pi^{-1}(D_i) = \{[x_0:x_1:x_2] \in \bar Y \mid x_i = 0\} = \{[0:1:\zeta^k], \, k = 0, \dots \tilde d -1\},
\]
where $\zeta$ is a primitive $2\tilde d$\=/ th unity root. Apply the map $\pi$ to these points. We~obtain
\[
D_0 = \{[0:1:\zeta\eta_0^k] , \, k = 0, \dots, dd_{12} - 1\},
\]
where $\eta_0$ is a primitive $dd_{12}$\=/ th unity root. One can get similarly
\begin{gather*}
D_1 = \{[1:0:\zeta\eta_1^k] , \, k = 0, \dots, dd_{02} - 1\},\\
D_2 = \{[0:\zeta\eta_2^k:1] , \, k = 0, \dots, dd_{01} - 1\},
\end{gather*}
where $\eta_1$ is a primitive $dd_{02}$\=/ th unity root and $\eta_2$ is a primitive $dd_{01}$\=/ th unity root.
\end{remark}
\section{Calculating examples}
We find proper polyhedral divisors for some specific classes of trinomial hypersurfaces in this section. Let us begin with rational hypersurfaces. Note that a trinomial hypersurface is rational if and only if the curve $Y$ from Theorem~\ref{main} is rational. This fact can be obtained from~\cite[Theorem~2.8]{VP}. In its turn, a curve is rational if and only if the genus of the curve is equal to zero. A genus of the curve from Theorem~\ref{main} is given by the following proposition.
\begin{proposition}\label{genus}\cite[Theorem~5.3]{ABHW}
Let $Y$ be the curve in $\P(d_{12}, d_{02}, d_{01})$ given by the equation
\[
w_0^{\tilde d / d_{12}} + w_1^{\tilde d / d_{02}} + w_2^{\tilde d / d_{01}} = 0.
\]
The genus of this curve is equal to
\[
g = \frac{d}{2}\left(\tilde d - (d_{01} + d_{02} + d_{12})\right) + 1.
\]
\end{proposition}
\begin{corollary}\cite[Lemma~5.6]{ABHW}
Let $Y$ be the curve from Theorem~\ref{main}. Then $Y$ is rational if and only if it belongs to one of the following types:
\begin{enumerate}[Type I.] 
\item
$d = 1$, $d_{01} = d_{02} = 1$, $d_{12} = s$
(up to renumbering),
\smallskip
\item
$d = 2$, $d_{01} = d_{02} = d_{12} = 1$.
\end{enumerate}
\end{corollary}
\subsection{Type I: Factorial case}
A trinomial hypersurface $X$ is factorial if and only if $d_0$, $d_1$, $d_2$ are pairwise coprime~\cite[Theorem~1.1(ii)]{HH}. According to Theorem~\ref{main}, the curve $Y$ where exists a $\sigma$\=/ polyhedral for a complexity one torus action on $X$ lies in $\P^2$ and is given by the equation
$w_0 + w_1 + w_2 = 0$.
Thus we obtain the following proposition for a factorial case.
\begin{proposition}
Let $X$ be a factorial trinomial hypersurface. Then the $\sigma$\=/ polyhedral divisor $\D$ for the complexity one torus action on $X$ has the form
\begin{gather*}
\D = \Delta_0 \cdot \{0\} + \Delta_1 \cdot \{1\} + \Delta_2 \cdot \{\infty\} \in \pp_\Q(\P^1, \sigma),\\
\Delta_0 = \conv\left(\left\{\left.S\left(\frac{1}{l_{0j}}e_k\right) \right| j = k, \,  k = 1, \dots , n_0\right\}\right) + \sigma,\\
\Delta_1 = \conv\left(\left\{\left.S\left(\frac{1}{l_{1j}}e_k\right) \right| j = k - n_0, \, k = n_0 + 1,\, \dots, \,n_0 + n_1\right\}\right) + \sigma,\\
\Delta_2 = \conv\left(\left\{\left.S\left(\frac{1}{l_{2j}}e_k\right) \right| j =  k - n_0 - n_1, \, k = n_0 + n_1 + 1, \,\dots,\, n\right\}\right) + \sigma,
\end{gather*}
where the cone $\sigma$ is defined like in Theorem~\ref{main}.
\end{proposition}
This result has been obtained earlier by Arzhantsev~\cite{Arzh}.
\begin{example}
Consider the hypersurface $X$ from Example~\ref{start}. Recall that a \mbox{$T$\=/ action} on $X$ is defined as
\[
(t_1, t_2) \cdot (T_{01}, T_{11}, T_{21}, T_{22}) = (t_1^5T_{01}, t_1^3T_{11}, t_2T_{21}, t_1^{15}t_2^{-1}T_{22}).
\]
Hence, matrices $F$ and $S$ look like
\[
F = \begin{pmatrix}5 & 3 & 0 & 15\\ 0 & 0& 1 & -1\end{pmatrix}^T, S = \begin{pmatrix}2 & -3 & 0 & 0\\0 & 0 & 1 & 0\end{pmatrix}.
\]
As a result, the divisor $\D$ has the form
\begin{multline*}
\D = \left(\left(\frac{2}{3}, 0\right) + \sigma \right) \cdot\{0\} + \left(\left(-\frac{3}{5}, 0\right) + \sigma \right) \cdot \{1\} + \\
+ (\{0\} \times [0, 1] + \sigma) \cdot \{\infty\} \in \pp_\Q(\P^1, \sigma),
\end{multline*}
where $\sigma = S(\Q^4_{\ge 0} \cap F(\Q^2)) = \c((1, 0), (1, 15))$.
\begin{figure}[ht]
\centering
\begin{tikzpicture}
\draw[black, thick] (-5.5, 0) -- (6.0, 0)
node[anchor=north] {$\P^1$};
\filldraw[black] (-5, 0) circle (2pt)
node[anchor=north] {0};
\filldraw[black!10] (-4.78, 1) -- (-4.65, 3) -- (-3.2, 2.8) -- (-3.1, 1) -- cycle;
\foreach \x in {-1, 0,..., 5}{\foreach \y in {-1, 0,..., 5}{\filldraw[black!20] (0.33 * \x - 5, 1 + 0.33 * \y) circle(0.5pt);}};
\draw[black, ->] (-5, 0.5) -- (-5, 3);
\draw[black, ->] (-5.5, 1) -- (-3, 1);
\draw[black] (-4.78, 1) -- (-4.65, 3) ++(0.5, -1.5) node {$\Delta_0$};
\filldraw[black] (0, 0) circle (2pt)
node[anchor=north] {1};
\filldraw[black!10] (-0.198, 1) -- (-0.065, 3) -- (1.8, 2.8) -- (1.9, 1) -- cycle;
\foreach \x in {-1, 0,..., 5}{\foreach \y in {-1, 0,..., 5}{\filldraw[black!20] (0.33 * \x, 1 + 0.33 * \y) circle(0.5pt);}};
\draw[black, ->] (0, 0.5) -- (0, 3);
\draw[black, ->] (-0.5, 1) -- (2, 1);
\draw[black] (-0.198, 1) -- (-0.065, 3) ++(0.5, -1.5) node {$\Delta_1$};
\filldraw[black] (5, 0) circle (2pt)
node[anchor=north] {$\infty$};
\filldraw[black!10] (5, 1) -- (5, 1.33) -- (5.11, 3) -- (6.8, 2.8) -- (6.9, 1) -- cycle;
\foreach \x in {-1, 0,..., 5}{\foreach \y in {-1, 0,..., 5}{\filldraw[black!20] (0.33 * \x + 5, 1 + 0.33 * \y) circle(0.5pt);}};
\draw[black, ->] (5, 0.5) -- (5, 3);
\draw[black, ->] (4.5, 1) -- (7, 1);
\draw[black] (5, 1.33) -- (5.11, 3) ++(0.5, -1.5) node {$\Delta_2$};
\end{tikzpicture}
\caption{}
\end{figure}
\end{example}
\subsection{Type I: Non-factorial case}
We can assume $d = 1$, $d_{01} = d_{02} = 1$, $d_{12} = s  \ge 2$ without loss of generality. According to Theorem~\ref{main}, the curve $Y$ where exists a $\sigma$\=/ polyhedral divisor for the complexity one torus action on $X$ lies in  $\P(s, 1, 1)$ and is given by the equation $w_0 + w_1^s + w_2^s = 0$. This curve can be associated with $\P^1$ by the following isomorphism:
\begin{gather*}
\phi_{\rn 1}: Y \ra \P^1,\, [w_0:w_1:w_2] \mapsto [w_1:w_2],\\
\phi_{\rn 1}^{-1}: \P^1 \ra Y, \, [z_0:z_1] \mapsto [-z_0^s - z_1^s:z_0:z_1].
\end{gather*}
Theorem~\ref{main} takes the following form.
\begin{proposition}
Let $X$ be a rational trinomial hypersurface of type I. Then the \mbox{$\sigma$\=/ polyhedral} divisor $\D$ for the complexity one torus action on $X$ has the form
\begin{gather*}
\D = \Delta_0 \cdot \sum_{k =0}^{s-1}\{\zeta^k\} + \Delta_1 \cdot \{0\} + \Delta_2 \cdot \{\infty\} \in \pp_\Q(\P^1, \sigma),\\
\Delta_0 = \conv\left(\left\{\left.S\left(\frac{1}{l_{0j}}e_k\right) \right| j = k, \,  k = 1, \dots , n_0\right\}\right) + \sigma,\\
\Delta_1 = \conv\left(\left\{\left.S\left(\frac{s}{l_{1j}}e_k\right) \right| j = k - n_0, \, k = n_0 + 1,\, \dots, \,n_0 + n_1\right\}\right) + \sigma,\\
\Delta_2 = \conv\left(\left\{\left.S\left(\frac{s}{l_{2j}}e_k\right) \right| j =  k - n_0 - n_1, \, k = n_0 + n_1 + 1, \,\dots,\, n\right\}\right) + \sigma,
\end{gather*}
where $\zeta$ is a primitive $s$\=/ th unity root and the cone $\sigma$ is defined like in Theorem~\ref{main}
\end{proposition}
\begin{example}
Consider a trinomial hypersurface $X$ of type I given by the equation $T_{01}^2 + T_{11}^3 + T_{21}^3T_{22}^3 = 0$. Matrices $L$, $F$ and $S$ have the form
\[
L = \begin{pmatrix}-2 & 3 & 0 & 0\\-2 & 0 & 3 & 3\end{pmatrix},\,
F = \begin{pmatrix}3 & 2 & 0 & 2\\0 & 0 & 1 & -1\end{pmatrix}^T,\,
S = \begin{pmatrix}1 & -1 & 0 & 0\\ 0 & 0 & 1 & 0\end{pmatrix}.
\]
As a result, the divisor $\D$ has the form
\begin{multline*}
\D = \left(\left(\frac{1}{2}, 0\right) + \sigma \right) \cdot (\{1\} + \{\zeta\} + \{\zeta^2\} ) + \left(\left(-\frac{1}{3}, 0\right) + \sigma \right) \cdot \{0\} +\\
+ \left(\{0\} \times \left[0, \frac{1}{3}\right] + \sigma\right) \cdot \{\infty\} \in \pp_\Q(\P^1, \sigma),
\end{multline*}
where $\zeta$ is a primitive $3$\=/ rd unity root and $\sigma = \c((1,0), (1,2))$.
\begin{figure}[ht]
\centering
\begin{tikzpicture}
\draw[black, thick] (-5.5, 0) -- (6.0, 0)
node[anchor=north] {$\P^1$};
\filldraw[black] (-5, 0) circle (2pt)
node[anchor=north] {1};
\filldraw[black] (-4.5, 0) circle (2pt)
node[anchor=north] {$\zeta$};
\filldraw[black] (-4, 0) circle (2pt)
node[anchor=north] {$\zeta^2$};
\filldraw[black!10] (-4.84, 1) -- (-3.84, 3) -- (-3.3 ,2.95) -- (-3.1, 1) -- cycle;
\foreach \x in {-1, 0,..., 5}{\foreach \y in {-1, 0,..., 5}{\filldraw[black!20] (0.33 * \x - 5, 1 + 0.33 * \y) circle(0.5pt);}};
\draw[black, ->] (-5, 0.5) -- (-5, 3);
\draw[black, ->] (-5.5, 1) -- (-3, 1);
\draw[black] (-4.84, 1) -- (-3.84, 3) ++(0.0, -1.5) node {$\Delta_0$};
\filldraw[black] (0, 0) circle (2pt)
node[anchor=north] {0};
\filldraw[black!10] (-0.11, 1) -- (0.89, 3) --(1.8, 2.9) -- (1.9, 1) -- cycle;
\foreach \x in {-1, 0,..., 5}{\foreach \y in {-1, 0,..., 5}{\filldraw[black!20] (0.33 * \x, 1 + 0.33 * \y) circle(0.5pt);}};
\draw[black, ->] (0, 0.5) -- (0, 3);
\draw[black, ->] (-0.5, 1) -- (2, 1);
\draw[black] (-0.11, 1) -- (0.89, 3) ++(0.0, -1.5) node {$\Delta_1$};
\filldraw[black] (5, 0) circle (2pt)
node[anchor=north] {$\infty$};
\filldraw[black!10] (5, 1) -- (5, 1.11) -- (5.95, 3) -- (6.8 ,2.95) -- (6.9, 1) -- cycle;
\foreach \x in {-1, 0,..., 5}{\foreach \y in {-1, 0,..., 5}{\filldraw[black!20] (0.33 * \x + 5, 1 + 0.33 * \y) circle(0.5pt);}};
\draw[black, ->] (5, 0.5) -- (5, 3);
\draw[black, ->] (4.5, 1) -- (7, 1);
\draw[black] (5, 1.11) -- (5.95, 3) ++(0.0, -1.5) node {$\Delta_2$};
\end{tikzpicture}
\caption{}
\end{figure}
\end{example}
\subsection{Rational hypersurfaces of type II}
In this case we have ${d_{01} = d_{02} = d_{12} = 1}$, ${d = 2}$. According to Theorem~\ref{main}, the curve $Y$ where exists a $\sigma$\=/ polyhedral divisor for the complexity one torus action on $X$ lies in $\P^2$ and given by the equation $w_0^2 + w_1^2 + w_2^2 = 0$. The genus of this curve is equal to zero hence the curve is isomorphic to $\P^1$. The isomorphism can be chosen as follows:
\begin{gather*}
\phi_{\rn 2}:Y \ra \P^1, \,[w_0:w_1:w_2] \mapsto 
\begin{cases}
[w_0 + \imath w_1:w_2],\, w_0 \not= -\imath w_1\\
[w_2:\imath w_1 - w_0], \, w_0 \not= \imath w_1
\end{cases},\\
\phi_{\rn 2}^{-1}: \P^1 \ra Y, \,[z_0: z_1] \mapsto [z_0^2 - z_1^2:-\imath (z_0^2+z_1^2):2z_0z_1],
\end{gather*}
where $\imath^2=-1$. Let us apply Theorem~\ref{main} to this case.
\begin{proposition}
Let $X$ be a rational trinomial hypersurface of type II. Then the \mbox{$\sigma$\=/ polyhedral} divisor $\D$ for the complexity one torus action on $X$ has the form
\begin{gather*}
\D = \Delta_0 \cdot (\{1\} + \{-1\}) + \Delta_1 \cdot (\{\imath\} + \{-\imath\}) + \Delta_2 \cdot (\{0\} + \{\infty\}) \in \pp_\Q(\P^1, \sigma),\\
\Delta_0 = \conv\left(\left\{\left.S\left(\frac{2}{l_{0j}}e_k\right) \right| j = k, \,  k = 1, \dots , n_0\right\}\right) + \sigma,\\
\Delta_1 = \conv\left(\left\{\left.S\left(\frac{2}{l_{1j}}e_k\right) \right| j = k - n_0, \, k = n_0 + 1,\, \dots, \,n_0 + n_1\right\}\right) + \sigma,\\
\Delta_2 = \conv\left(\left\{\left.S\left(\frac{2}{l_{2j}}e_k\right) \right| j =  k - n_0 - n_1, \, k = n_0 + n_1 + 1, \,\dots,\, n\right\}\right) + \sigma,
\end{gather*}
where $\imath^2=-1$ and the cone $\sigma$ is defined as in Theorem~\ref{main}. 
\end{proposition}
\begin{example}
Consider a trinomial hypersurface $X$ of type I given by the equation $T_{01}^2 + T_{11}^4 + T_{21}^2T_{22}^4 = 0$. Matrices $L$, $F$ and $S$ have the form
\[
L = \begin{pmatrix}-2 & 4 & 0 & 0\\-2 & 0 & 2 & 4\end{pmatrix},\,
F = \begin{pmatrix}2 & 1 & 0 & 1\\0 & 0 & -2 & 1\end{pmatrix}^T,\,
S = \begin{pmatrix}1 & -1 & 0 & 0\\ 0 & -1 & 0 & 1\end{pmatrix}.
\]
As a result, the divisor $\D$ has the form
\begin{multline*}
\D = \left((1, 0)+ \sigma\right)\cdot (\{1\} + \{-1\}) + \left(\left(-\frac{1}{2}, -\frac{1}{2}\right)+ \sigma\right)\cdot (\{\imath\} + \{-\imath\})+ \\ 
+ \left(\{0\} \times \left[0, \frac{1}{2}\right] + \sigma\right)\cdot (\{0\} + \{\infty\}) \in \pp_\Q(\P^1, \sigma),
\end{multline*}
where $\imath^2=-1$ and $\sigma = \c((1, 0), (1, -1))$.
\begin{figure}[ht]
\centering
\begin{tikzpicture}
\draw[black, thick] (-5.5, 0) -- (6.6, 0)
node[anchor=north] {$\P^1$};
\filldraw[black] (-5, 0) circle (2pt)
node[anchor=north] {1};
\filldraw[black] (-4.4, 0) circle (2pt)
node[anchor=north] {$-1$};
\filldraw[black!10] (-4.66, 2.5) -- (-3.16, 1) -- (-3.1, 2.5) -- cycle;
\foreach \x in {-1, 0,..., 5}{\foreach \y in {1, 0,..., -5}{\filldraw[black!20] (0.33 * \x - 5, 2.5 + 0.33 * \y) circle(0.5pt);}};
\draw[black, ->] (-5, 0.5) -- (-5, 3);
\draw[black, ->] (-5.5, 2.5) -- (-3, 2.5);
\draw[black] (-4.66, 2.5) -- (-3.16, 1) ++(-0.33, 1) node {$\Delta_0$};
\filldraw[black] (0, 0) circle (2pt)
node[below=2pt] {$\imath$};
\filldraw[black] (0.6, 0) circle (2pt)
node[below=1pt] {$-\imath$};
\filldraw[black!10] (-0.16, 2.34) -- (1.63, 0.55) -- (1.9, 2.34) -- cycle;
\foreach \x in {-1, 0,..., 5}{\foreach \y in {1, 0,..., -5}{\filldraw[black!20] (0.33 * \x, 2.5 + 0.33 * \y) circle(0.5pt);}};
\draw[black, ->] (0, 0.5) -- (0, 3);
\draw[black, ->] (-0.5, 2.5) -- (2, 2.5);
\draw[black] (-0.16, 2.34) -- (1.63, 0.55) ++(-0.33, 1) node {$\Delta_1$};
\draw[black] (-0.16, 2.34) -- (1.9, 2.34);
\filldraw[black] (5, 0) circle (2pt)
node[below=1pt] {0};
\filldraw[black] (5.6, 0) circle (2pt)
node[below=2pt] {$\infty$};
\filldraw[black!10] (5, 2.5) -- (5, 2.17) -- (6.62, 0.55) -- (6.9, 2.5) -- cycle;
\foreach \x in {-1, 0,..., 5}{\foreach \y in {1, 0,..., -5}{\filldraw[black!20] (0.33 * \x + 5, 2.5 + 0.33 * \y) circle(0.5pt);}};
\draw[black, ->] (5, 0.5) -- (5, 3);
\draw[black, ->] (4.5, 2.5) -- (7, 2.5);
\draw[black] (5, 2.17) -- (6.62, 0.55) ++(-0.33, 1) node {$\Delta_2$};
\end{tikzpicture}
\caption{}
\end{figure}
\end{example}
\subsection{Pham-Brieskorn surfaces}
Let $n_0 = n_1 = n_2 = 1$. Then $X$ is given by the equation $x_0^{d_0} + x_1^{d_1} + x_2^{d_2} = 0$ in $\mathbb{A}^3$. Such varieties are well known as Pham-Brieskorn surfaces. There is an action of a one-dimensional torus $T$ of complexity one on $X$. One can choose $F$ and $S$ so that $\sigma = \Q_{\ge 0}$. For this it is enough to make all elements of $F$ non-negative. The map $F$ can be described by the vector
\[
f = (\bar d/d_0 , \bar d/d_1 , \bar d/d_2),
\]
where $\bar d := \lcm(d_0, d_1, d_2)$. Coordinates of $f$ are coprime. Hence, there exists such a vector $s = (s_0 , s_1 , s_2)$ that $(f, s) = 1$. Thus, the map $S$ can be described by $s$.
\begin{proposition}
Let $X$ be a Pham-Brieskorn surface given by the equation
${x_0^{d_0} + x_1^{d_1} + x_2^{d_2} = 0}$.
Then the \mbox{$\sigma$\=/ polyhedral} divisor $\D$ for a complexity one torus action on $X$ has the form
\begin{gather*}
\D = \Delta_0 \cdot D_0 + \Delta_1 \cdot D_1 + \Delta_2 \cdot D_2 \in \pp_\Q(Y, \Q_{\ge 0}),\\
\Delta_0 = \left(s_0\frac{dd_{01}d_{02}}{d_0}, +\infty\right), \, \Delta_1 = \left(s_1\frac{dd_{01}d_{12}}{d_1}, +\infty\right), \,
\Delta_2 = \left(s_2\frac{dd_{02}d_{12}}{d_2}, +\infty\right),
\end{gather*}
where the cone $\sigma$ and the curve $Y$ are defined as in Theorem~\ref{main}. 
\end{proposition}
\begin{example}
Consider a Pham-Brieskorn surface $X$ given by the equation $x_0^2 + x_1^3 + x_2^6 = 0$ in $\K^3$.
Vectors $f$ and $s$ looks like
\[
f = (3, \,2, \,1),\, s = (1, \,-1, \,0).
\]
The curve $Y$ lies in $P(3, 2, 1)$ and the genus of $Y$ is equal to one, according to Proposition~\ref{genus}. Hence, we find out that the divisor $\D$ has the form
\begin{gather*}
\D = \left[1, +\infty\right) \cdot D_0+ [-1, +\infty) \cdot D_1,\\
D_0 = \sum_{k = 1}^3\{[0:\zeta^k:1]\}, \, D_1 = \{[\imath:0:1]\} + \{[-\imath:0:1]\} \in \pp_\Q(Y, \Q_{\ge 0}),
\end{gather*}
where $\zeta$ is a primitive 3\=/ rd unity root.
\end{example}
\begin{remark}
This example shows us that some of polyhedral coefficients of $\D$ can be equal to $\sigma$.
\end{remark}
\subsection{Non-rational case}
In conclusion we find a pp-divisor for a non-rational hypersurface.
\begin{example}
Consider a non-rational trinomial hypersurface $X$ given by the equation $T_{01}^2 + T_{11}^4 + T_{21}^2T_{22}^4 = 0$. Matrices $L$, $F$ and $S$ have the form
\[
L = \begin{pmatrix}-2 & 3 & 0 & 0 \\ -2 & 0 & 6 & 6\end{pmatrix}, \,
F = \begin{pmatrix}3 & 2 & 0 & 1\\ 0 & 0 & 1 & -1\end{pmatrix}^T,\, 
S = \begin{pmatrix}1 & -1 & 0 & 0\\ 0 & 0 & 1 & 0\end{pmatrix}. 
\]
\tikzset{
    partial ellipse/.style args={#1:#2:#3}{
        insert path={+ (#1:#3) arc (#1:#2:#3)}
    }
}
\begin{figure}[h]
\centering
\begin{tikzpicture}
\draw[thick] (0, 1) ellipse (3 and 1.5);
\draw[thick] (0, 2) [partial ellipse=210:330:2.5 and 1.5];
\draw[thick] (0, 0) [partial ellipse=42.4:137.6:2.5 and 1.5];
\filldraw (0, 2) circle (2pt) node[anchor=north] {$D_0$};
\filldraw (0.4, 2.2) circle (2pt);
\filldraw (-0.4, 2.2) circle (2pt);
\filldraw (-1.6, 0.385) circle (2pt)  node[anchor=east] {$D_2$};
\filldraw (1.85, 0.585) circle (2pt) ;
\filldraw (2.25, 0.79) circle (2pt) ++ (-0.08, -0.15) node[anchor=north] {$D_1$};

\filldraw[black!10] (3.17, 0) -- (5.14, 1.97) -- (5.4, 0) -- cycle;
\foreach \x in {-1, 0,..., 5}{\foreach \y in {-1, 0,..., 5}{\filldraw[black!20] (0.33 * \x + 3.5, 0.33 * \y) circle (0.5pt);}};
\draw[->] (3, 0) -- (5.5, 0);
\draw[->] (3.5, -0.5) -- (3.5, 2);
\draw (3.17, 0) -- (5.14, 1.97) ++ (-0.5, -1.23) node {$\Delta_1$};
\filldraw[black!10] (-5, 0) -- (-5, 0.33) -- (-3.33, 2.0) -- (-3.1, 0) -- cycle;
\foreach \x in {-1, 0,..., 5}{\foreach \y in {-1, 0,..., 5}{\filldraw[black!20] (0.33 * \x - 5.0, 0.33 * \y) circle (0.5pt);}};
\draw[->] (-5.5, 0) -- (-3, 0);
\draw[->] (-5, -0.5) -- (-5, 2);
\draw (-5, 0.33) -- (-3.33, 2.0) ++(-0.53, -1.23) node {$\Delta_2$};
\filldraw[black!10] (-0.17, 3.5) -- (1.35, 5.02) -- (1.4, 3.5) -- cycle;
\foreach \x in {-1, 0,..., 5}{\foreach \y in {-1, 0,..., 5}{\filldraw[black!20] (0.33 * \x - 0.5, 3.5 + 0.33 * \y) circle(0.5pt);}};
\draw[->] (-1, 3.5) -- (1.5, 3.5);
\draw[->] (-0.5, 3) -- (-0.5, 5.5);
\draw (-0.17, 3.5) -- (1.35, 5.02) ++ (-0.5, -1.1) node {$\Delta_0$};

\end{tikzpicture}
\caption{}
\end{figure}
\newpage
The curve $Y$ is similar to the previous example. In this case the $\sigma$\=/ polyhedral divisor has the~form
\begin{gather*}
\D = \left((1, 0) + \sigma\right) \cdot D_0+ ((-1, 0) + \sigma) \cdot D_1 + (\{0\} \times [0,1] + \sigma) \cdot D_2,\\
D_0 = \sum_{k = 1}^3\{[0:\zeta^k:1]\}, \, D_1 = \{[\imath:0:1]\} + \{[-\imath:0:1]\},\, D_2 = \{[1:-1:0]\},
\end{gather*}
where
$\sigma = \c \left( (1, 0), \, (1, 1)\right)$, 
$\imath^2 = -1$, and $\zeta$ is a primitive 3\hyph rd unity root.

\end{example}

\bibliographystyle{amsplain}
\bibliography{bibliography_short_titles}

\providecommand{\bysame}{\leavevmode\hbox to3em{\hrulefill}\thinspace}
\providecommand{\MR}{\relax\ifhmode\unskip\space\fi MR }
\providecommand{\MRhref}[2]{%
  \href{http://www.ams.org/mathscinet-getitem?mr=#1}{#2}
}
\providecommand{\href}[2]{#2}
\begin{thebibliography}{10}

\bibitem{AH}
Klaus Altmann and J\"{u}rgen Hausen, \emph{Polyhedral divisors and algebraic
  torus actions}, Math. Ann. \textbf{334} (2006), no.~3, 557--607.

\bibitem{TVar}
Klaus Altmann, Nathan~Owen Ilten, Lars Petersen, Hendrik S\"{u}\ss, and Robert
  Vollmert, \emph{The geometry of {$T$}-varieties}, Contributions to Algebraic
  Geometry, EMS Ser. Congr. Rep., Eur. Math. Soc., Z\"{u}rich, 2012, 17--69.

\bibitem{Arzh}
Ivan Arzhantsev, \emph{On rigidity of factorial trinomial hypersurfaces}, Int.
  J. Algebr. Comput. \textbf{26} (2016), no.~5, 1061--1070.

\bibitem{ABHW}
Ivan Arzhantsev, Lukas Braun, J\"{u}rgen Hausen, and Milena Wrobel, \emph{Log
  terminal singularities, platonic tuples and iteration of {C}ox rings}, Eur.
  J. Math. \textbf{4} (2018), no.~1, 242--312.

\bibitem{CLS}
David Cox, John Little, and Henry Schenck, \emph{Toric varieties}, Grad. Stud.
  Math., vol. 124, Amer. Math. Soc., Providence, RI, 2011.

\bibitem{F}
William Fulton, \emph{Introduction to toric varieties}, Ann. of Math. Studies,
  vol. 131, Princeton University Press, Princeton, NJ, 1993, The William H.
  Roever Lectures in Geometry.

\bibitem{HH}
J\"urgen Hausen and Elaine Herppich, \emph{Factorially graded rings of
  complexity one}, Torsors, \'etale homotopy and applications to rational
  points, London Math. Soc. Lecture Note Ser., vol. 405, Cambridge Univ. Press,
  Cambridge, 2013, 414--428.

\bibitem{VP}
Vladimir Popov and \`{E}rnest Vinberg, \emph{Invariant theory}, Encycl. of
  Math. Sci., 123--278, Springer Berlin Heidelberg, Berlin, Heidelberg, 1994.

\bibitem{WPS}
Michele Rossi and Lea Terracini, \emph{Linear algebra and toric data of
  weighted projective spaces}, Rend. Semin. Mat. Univ. Politec. Torino
  \textbf{70} (2012), no.~4, 469--495.

\bibitem{Vin}
\`{E}rnest Vinberg, \emph{Complexity of action of reductive groups}, Funct.
  Anal. Appl. \textbf{20} (1986), no.~1, 1--11.

\end{thebibliography}
\end{document}